%% file: article.tex
\documentclass[finalversion]{FPSAC2024}
\articlenumber{65}

\usepackage{tikz-qtree}

\newtheorem{thm}{Theorem}
\newtheorem{prop}[thm]{Proposition}
\newtheorem{lem}[thm]{Lemma}

\theoremstyle{definition}

\newtheorem{defi}[thm]{Definition}

\newtheorem{exa}[thm]{Example}

\usepackage{tikz}
\usepackage{tikz-cd}
\usepackage{wrapfig}
\usepackage{tikzit}
\input{article.tikzstyles}
\input{article.tikzdefs}
\usepackage{pgfplots}
\pgfplotsset{compat=1.13}
\usepackage{caption,subcaption}
\usepackage{mathtools}
\usepackage[normalem]{ulem}

\def\Z{\mathbb{Z}}
\def\N{\mathbb{N}}
\def\C{\mathbb{C}}
\def\K{\mathbb{K}}
\def\Q{\mathbb{Q}}
\def\Cmul{\mathbb{C}_{\mathrm{mul}}}
\def\calS{\mathcal{S}}
\def\invS{1/S(x,y)}
\def\calW{\mathcal{W}}
\def\calH{\mathcal{H}}
\def\Gmod{\mathcal{G}_{\lambda}}
\def\calG{\mathcal{G}}
\def\Ktld{\widetilde{K}}
\def\cO{\mathcal{O}}
\def\id{\operatorname{id}}

\def\Aut{\operatorname{Aut}}
\def\kinv{k_{\mathrm{inv}}}

\title[A Galois structure on the orbit of large steps walks in the quadrant]{A Galois structure on the orbit \\ of large steps walks in the quadrant}
\author{Pierre Bonnet\thanks{\href{mailto:pierre.bonnet@u-bordeaux.fr}{pierre.bonnet@u-bordeaux.fr}}\addressmark{1}
\and Charlotte Hardouin\thanks{\href{mailto:hardouin@math.univ-toulouse.fr}{hardouin@math.univ-toulouse.fr}}\addressmark{2}}

\address{\addressmark{1}LaBRI, Université de Bordeaux, Bordeaux \\ \addressmark{2}Institut de mathématiques, Université Paul Sabatier, Toulouse}

\received{\today}

\revised{}

\abstract{The enumeration of weighted walks in the quarter plane reduces to studying
  a functional equation with two catalytic variables. When the steps of the walk are small,
  Bousquet-Mélou and Mishna defined a group called \emph{the group of
    the walk} which turned out to be crucial in the classification of the small steps models.
  In particular, its action on the catalytic variables provides a convenient set of changes
  of variables in the functional equation. This particular set
  called the \emph{orbit} has been generalized to models with arbitrary large steps
  by Bostan, Bousquet-Mélou and Melczer (BBMM).
  However, the orbit had till now no underlying group.

  In this article, we endow the orbit with the action of a Galois group, which extends
  the notion of the group of the walk to models with large steps. As an application, we
  look into a general strategy to prove the algebraicity of models with small backwards
  steps, which uses the
  fundamental objects that are \emph{invariants} and \emph{decoupling}. The group action
  on the orbit allows us to develop a Galoisian approach to these two notions. Up
  to the knowledge of the finiteness of the orbit, this gives
  systematic procedures to test their existence and construct them. Our constructions
  lead to the first proofs of algebraicity of weighted models with large steps,
  proving in particular a conjecture of BBMM, and allowing to find
  new algebraic models with large steps.
}


\usepackage{biblatex}
\begin{document}

\maketitle

\section{Introduction and preliminaries} \label{sec:intr-prel}

\subsection{Walks in the quarter plane\label{sec:walks-quarter-plane}}
A weighted walk in the quarter plane is defined as follows.
Consider a finite subset $\calS$ of $\Z\times\Z$. To each \emph{step} $s$
of $\calS$ we attach its \emph{weight} $w_s$ which is a nonzero complex number.
The tuple $\calW=(\calS,(w_s)_{s \in \calS})$ is called a
\emph{weighted model of walks}.
A \emph{weighted walk in the quarter plane} of length $n$
on the model $\calW$ is then a sequence of points $P_0, \dots, P_n$ in $\N
\times \N$ such that for all $i$ there exists $s_i \in \calS$
satisfying $s_i + P_i = P_{i+1}$. The \emph{weight} of the walk
is the product $w_{s_0} w_{s_1} \dots w_{s_{n-1}}$ of the weights of the steps
taken by the walk.

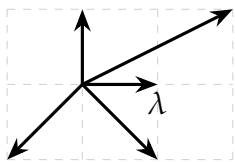
\begin{figure}[ht]
    \centering
    \begin{subfigure}[c]{0.6\textwidth}
      \centering
      \scalebox{1}{\input{g2_mod.tikz}}
      \caption*{The weighted model $\Gmod$
        (for which we take any nonzero $\lambda$ in $\C$)
          along with an example of a walk on $\Gmod$,
          of length $8$, ending at $(3,0)$ and of
        weight $ \lambda^2$. }
        \label{fig:g2_mod}
    \end{subfigure}
    \begin{subfigure}[c]{0.3\textwidth}
      \centering
      \scalebox{0.8}{
        \input{ex_walk_g2.tikz}
        }
        \label{fig:walk2D}
    \end{subfigure}
    \caption{An example of weighted model and walk} \label{fig:fig_gmod_walk}

\end{figure}
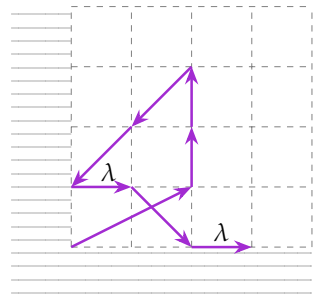

The enumeration of weighted walks in the quarter plane
has attracted a lot of attention over the past 20 years.
Indeed, these objects are general enough to encode many objects in
combinatorics (families of permutations, trees, maps),
probability theory (stochastic processes, games of chance, sums of
discrete random variables) or statistics (non-parametric tests).
The attraction for this topic comes from the fact that the solution
of this problem requires many different techniques and
points of view, from combinatorics of course, but also from probability
theory, computer algebra, differential Galois theory, complex analysis, geometry \ldots

\subsection{Generating function and classification}\label{sec:gener-funct-class-1} 
Given a model $\calW$, denote by $q_{i,j,n}$ the sum of the weights of walks in the
quadrant on $\calW$ of length $n$ starting at $P_0$ (taken as $(0,0)$ unless stated
otherwise) and terminating at $(i,j)$.
The generating function for these walks is defined as 
\[
  Q(X,Y,t) = \sum_{i,j,n \ge 0} q_{i,j,n} X^i Y^j t^n.
\]
The weighted model is encoded as the \emph{Laurent polynomial of the model}
  defined as $S(X,Y) = \sum_{s \in \calS} w_s X^{s_x} Y^{s_y}$.
From any weighted model $\calW$, it is quite easy to form a functional
equation for $Q(X,Y,t)$, as we demonstrate in the following example.
\begin{exa}\label{ex:func_eq_gmod}
  Let $\Gmod = \{(-1,-1), (0,1), ((1,0),\lambda), (2,1), (1,-1)\}$ as in
  Figure~\ref{fig:fig_gmod_walk} (Example~2.1 in \cite{BHorb}, see also Remark~2.2 for
  alternate weightings).
  Its Laurent polynomial is
  $S(X,Y) = \frac{1}{XY} + Y + \lambda X + \frac{X}{Y} + X^2 Y$.
  We now construct a recurrence on the walks:
  a walk terminating at coordinates $(i,j)$ can be completed by a step $s$ of $\calS$
  as long as $(i,j) + s$ is in $\N\times\N$. This translates into the following functional equation:
  \begin{align*}
    & Q(X,Y,t) = 1 + t X^2 Y Q(X,Y,t) + t \lambda X Q(X,Y,t) + t Y Q(X,Y,t) \\
             &+ t \frac{X}{Y} \left(Q(X,Y,t) - Q(X,0,t) \right)
             + t \frac{1}{XY} \left(Q(X,Y,t) - Q(X,0,t) - Q(0,Y,t) + Q(0,0,t) \right).
  \end{align*}
  Such an equation is then usually put in the following normal form:
\begin{equation}\label{eq:eqfunc_g2}
\Ktld(X,Y,t) Q(X,Y,t) = XY - t (X^2 + 1) Q(X,0,t) - tQ(0,Y,t) +
t Q(0,0,t),
\end{equation}
with $\Ktld(X,Y,t)$ the \emph{kernel polynomial of the walk}
being equal here to $XY \left(1 - tS(X,Y)\right)$.
\end{exa}

Given a class of combinatorial objects, a natural question is to
determine where its generating function fits in the classical hierarchy of power
series of $\C(X,Y)[[t]]$
\[\emph{rational} \subset \emph{algebraic} \subset \emph{D-finite} \subset \emph{D-algebraic}, \]
  where algebraic series satisfy polynomial equations; D-finite series
  satisfy one linear differential equation in each variable $X$, $Y$, $t$;
  and D-algebraic series satisfy polynomial differential equations, all the coefficients
  being taken in the polynomial ring $\C[X,Y,t]$.

For walks, this hierarchy measures the complexity of a
model: the lower its generating function in this hierarchy, the simpler the walks.
The \emph{catalytic variables equations} like (\ref{eq:eqfunc_g2})
do not immediately allow to conclude to the position of their solutions in this
hierarchy, hence the question of classifying the complexity of a model of walks
in the quarter plane is highly nontrivial.

\begin{figure}[ht]
  \begin{subfigure}[c]{0.5\linewidth}
    \centering
    \scalebox{1.2}{
      \input{clas_small.tikz}
    }
    \caption{$\calS \subset \{-1,0,1\}^2$} \label{fig:clas_small_steps}
  \end{subfigure}
  \begin{subfigure}[c]{0.4\linewidth}
    \centering
    \scalebox{1.2}{
    \input{clas_large_steps.tikz}
    }
    \caption{$\calS \subset \{-2,-1,0,1\}^2$} \label{fig:clas_large_steps}
  \end{subfigure}
  \caption{The partial classifications for two families of unweighted models}
\end{figure}
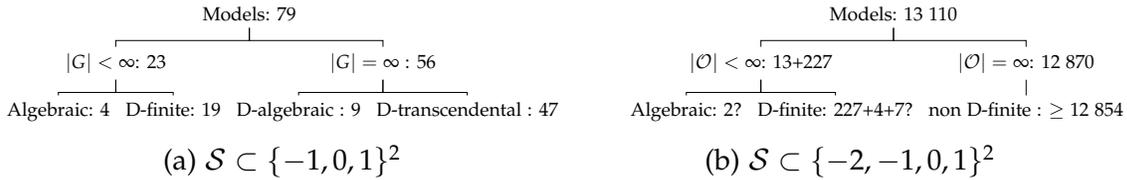

For instance, consider the restriction of the problem to unweighted models
with steps contained in $\{-1,0,1\}^2$ (these models are
commonly called \emph{small steps models}). The classification completed
in 2018 and summarized in Figure~\ref{fig:clas_small_steps} shows the numerous
different behaviours that arise.
Following the success of this first classification, the set of
considered models has been extended to models with steps in
$\{-2,1,0,-1\}^2$ in~\cite{bostan2018counting}. It is summarized in
Figure~\ref{fig:clas_large_steps}, and is currently incomplete. One of
the reasons is that not all the tools used in the classification of
small steps models extend to large steps.

For small steps models, \cite[Chapter 4]{FIM} proposes a complete strategy
with Galois theoretic tools to classify solutions of functional
equation of the form \eqref{eq:eqfunc_g2} when the kernel polynomial is biquadratic
and the orbit is finite. Unfortunately, theses tools rely heavily on
an elliptic uniformization of the algebraic curve associated with the
kernel polynomial. We propose here a
an extension of these tools with in sight a particular strategy to prove algebraicity,
which goes beyond the elliptic framework since it deals with kernel
polynomials of arbitrary degree. It has the advantage to stay within
the realm of Laurent power series and is almost algorithmic until an
algebraic characterization of certain invariants. Moreover, we hope
that the geometric framework hidden behind our constructions will
allow to adapt entirely the strategy of \cite{BBMR16} to large steps models.

\subsection{The group and the orbit}\label{sec:intro-group-orbit} 

A fundamental object which arises in the study of models with small steps is
the \emph{group of the walk}, introduced by Bousquet-Mélou and Mishna
in~\cite{BMM}, following \cite{FIM}. It is defined as follows. For a small
steps model,
we can write its Laurent polynomial in two ways:
\begin{align*}
  S(X,Y) &= A_{-1} (X)/ Y + A_{0}(X) + A_{1}(X) Y \\
         &= B_{-1}(Y)/ X + B_{0}(Y) + B_{1}(Y) X.
\end{align*}
Assume $A_{-1}(X)$, $A_{1}(X)$, $B_{-1}(Y)$ and $B_{1}(Y)$ to be nonzero.
The polynomial $S(x,y)$ is left unchanged by the
two \emph{birational transformations} (that are involutions) of $\C \times \C$
defined as \[
  \Phi \colon (u,v) \mapsto \left(\frac{B_{-1}(v)}{u B_{1}(v)},v\right) \qquad
  \Psi \colon (u,v) \mapsto \left(u,\frac{A_{-1}(u)}{v A_1(u)}\right). \]
The \emph{group of the walk} is then defined as $\left<\Phi,\Psi\right>$,
the subgroup of birational transformations of $\C \times \C$ generated by $\Phi$ and $\Psi$.
This group turned out to be a crucial algebraic invariant of a model with small steps.
For instance, an unweighted model with small steps
has a D-finite generating function if and only if the group is finite
(see the introduction of~\cite{BBMR16}).

The group also acts on pairs of catalytic variables. The orbit of its action on the
pair $(x,y)$ has a graph structure: the vertices are the pairs $(u,v)$ of the
orbit of $(x,y)$, and two pairs are adjacent if one can be obtained from the
other by applying $\Phi$ or $\Psi$ to it. 

If the model contains a large step, the equation $S(x,y)=S(x,y')$
may have non-rational solutions in $x$ and $y$ because of the higher
degree of the polynomials, therefore in that case this group cannot be defined
as a group of birational transformations of $\C \times \C$. Nonetheless,
Bostan, Bousquet-Mélou and Melczer noted in~\cite{bostan2018counting}
that the graph could be defined independently from the group,
and called it the \emph{orbit of the walk}. It is defined
as follows.
Denote by $\K = \overline{\C(x,y)}$ the algebraic closure of $\C(x,y)$
for two indeterminates $x$ and $y$.

\begin{defi}[Definition~3.1 in~\cite{bostan2018counting}]
  Let $(u,v)$ and $(u',v')$ be in $\K \times \K$. Then $(u,v)$ and $(u',v')$ are called
  \emph{$x$-adjacent} if $S(u,v) = S(u',v')$ and $u=u'$. Similarly, they are called
  \emph{$y$-adjacent} if $S(u,v) = S(u',v')$ and $v=v'$. We denote these two
  equivalence relations by $\sim^x$ and $\sim^y$. Two pairs are then adjacent
  if they are either $x$-adjacent or $y$-adjacent, and this relation is denoted by $\sim$.
  We denote by $\sim^*$ the transitive closure of $\sim$. The \emph{orbit of the walk} is
  the set of pairs $(u,v)$ such that $(u,v) \sim^* (x,y)$, and is denoted
  $\cO$.
\end{defi}

\begin{exa} \label{exa:orb_krew}
  For the small steps model $S(X,Y) = 1/X + 1/Y + XY$, the orbit
  is the cycle
  \[(x,y)
    \xleftrightarrow{\Phi} \left(\frac{1}{xy},y\right)
    \xleftrightarrow{\Psi} \left(\frac{1}{xy},x\right)
    \xleftrightarrow{\Phi} \left(y,x\right)
    \xleftrightarrow{\Psi} \left(y,\frac{1}{xy}\right)
    \xleftrightarrow{\Phi} \left(x,\frac{1}{xy}\right)
    \xleftrightarrow{\Psi} \left(x,y\right) \]
\end{exa}
\begin{exa} \label{exa:orb_o12}
  For $\Gmod$, the equation $S(x,y) = S(x',y)$ has three solutions $x$, $x_1$ and $x_2$
  with $x_1$ algebraic of degree two over $\C(x,y)$. Continuing the construction, it turns out
  that the orbit $\cO$ is finite of size $12$ with coordinates in $\C(x,y,x_1)$
  (see Figure~\ref{fig:orb_o12}).
\end{exa}

In Section~\ref{sec:galo-struct-orbit}, we introduce a Galois framework to study the orbit.
It allows us to generalize the group of the walk to weighted models with arbitrarily large steps.

\subsection{A strategy based on invariants and decoupling}\label{sec:invar-deco} 

In their article~\cite{BBMR16}, Bernardi, Bousquet-Mélou and Raschel introduced a general
strategy for proving algebraicity of a small steps model in the quadrant, later adapted
to models of the three-quadrant cone in~\cite{BousquetMelouThreequadrant}. This method
relies critically on objects called \emph{invariants}, of which we use one flavor
  in the realm of formal power series, the \emph{t-invariants}.

Define the subring $\Cmul(X,Y)[[t]]$ of power series of $\C(X,Y)[[t]]$ whose
coefficients in the $t$-expansion are of the form $\frac{A_n(X,Y)}{B_n(X)C_n(Y)}$
for $A_n$, $B_n$ and $C_n$ polynomials over~$\C$. A power series in $\Cmul(X,Y)[[t]]$
is said to have \emph{poles of bounded order at 0} if there is a bound on the order
of the poles at $X=0$ and $Y=0$ of its coefficients (Definition 2.1 in~\cite{BousquetMelouThreequadrant}).

\begin{defi} \label{defi:t_equiv}
  Let $F(X,Y,t)$ and $G(X,Y,t)$ be two power series of $\Cmul(X,Y)[[t]]$. They are
  \emph{$t$-equivalent} (with respect to $\Ktld$) if
  the power series $\frac{F(X,Y,t)-G(X,Y,t)}{\Ktld(X,Y,t)}
   \in \Cmul(X,Y)[[t]]$
  has poles of bounded order at~$0$. This equivalence is denoted by $F(X,Y,t) \equiv G(X,Y,t)$.
\end{defi}

The $t$-equivalence relation is compatible with ring operations,
and is used to define the notions of \emph{$t$-invariants} and \emph{$t$-decoupling}.

\begin{defi}[Invariants, Def.~2.3 in~\cite{BousquetMelouThreequadrant}] \label{def:invariants}
  Let $F(X,t), G(Y,t) \in \Cmul(X,Y)[[t]]$. The pair $(F(X,t),G(Y,t))$ is called a \emph{pair of $t$-invariants}
  if $F(X,t) \equiv G(Y,t)$.
\end{defi}

\begin{defi}[Decoupling] \label{def:decoupling}
  Let $H(X,Y,t)$ be a series of $\Cmul(X,Y)[[t]]$.
  Then $H(X,Y,t)$ admits a \emph{$t$-decoupling} if there exist $F(X,t)$ in $\Cmul(X)[[t]]$
  and $G(Y,t)$ in $\Cmul(Y)[[t]]$ such that $H(X,Y,t) \equiv F(X,t) + G(Y,t)$.
\end{defi}

\begin{exa}
  Consider the model defined by $S(X,Y) = XY + \frac{1}{X} + \frac{1}{Y}$
  (the same as Example~\ref{exa:orb_krew}).
  The fraction $XY$ admits the obvious decoupling $XY \equiv \frac{1}{t} - \frac{1}{X} - \frac{1}{Y}$.
  Moreover, the following identity induces a pair of rational invariants:
$X + \frac{1}{tX} - \frac{1}{X^2} \equiv Y + \frac{1}{tY} - \frac{1}{Y^2}$.
\end{exa}

When all the components of a pair of $t$-invariants or a $t$-decoupling are
  rational fractions, we speak of \emph{rational} $t$-invariants or $t$-decoupling.
  This notion of invariants intervenes in the following result, on which
  the strategy of \cite{BBMR16,BousquetMelouThreequadrant} relies crucially.
\begin{lem}[Lemma~2.6 in~\cite{BousquetMelouThreequadrant}] \label{lem:inv_lemma}
  Let $(F(X,t),G(Y,t))$ be a pair of $t$-invariants. If the coefficients
  of the power series $\frac{F(X,t)-G(Y,t)}{\Ktld(X,Y,t)} \in \Cmul(X,Y)[[t]]$
  have no pole at $X=0$ nor $Y=0$, then there exists a series $A(t)$ in $\C[[t]]$ such that
  $F(X,t) = G(Y,t) = A(t)$.
\end{lem}

The strategy of \cite{BBMR16,BousquetMelouThreequadrant} applies verbatim to large steps models
with small backward steps,
and goes as follows. Using a \emph{rational $t$-decoupling} of $XY$
(or more generally $X^{k+1} Y^{l+1}$ for other starting points) and
the special shape of the equation (e.g. Equation~(\ref{eq:eqfunc_g2})),
we construct a first pair of $t$-invariants. Next, we combine it with a
pair of \emph{non-constant rational $t$-invariants} using ring operations,
to eventually obtain a third pair of invariants that satisfy the conditions
of Lemma~\ref{lem:inv_lemma}. As this pair involves $Q(X,0,t)$ and $Q(0,Y,t)$,
the lemma gives two equations with one catalytic variable on these series.
By a result of Bousquet-Mélou and Jehanne in~\cite{BMJ}, they must be algebraic,
so is $Q(X,Y,t)$.

The existence of non-constant rational $t$-invariants and decoupling
is crucial to conduct this strategy. In sections \ref{sec:galois-invariants} and \ref{sec:decoupling},
we give a Galois approach to these two notions, which exploits
the notion of the group of the walk introduced in
Section~\ref{sec:galo-struct-orbit}, providing a \emph{systematic}
construction of these objects up to their existence and the finiteness of the orbit.
This is an alternative approach to the one developed in \cite{BuchacherKauersPogudin}
  where the authors search for a polynomial decoupling (see the discussion of
  Example~3.19 in \cite{BHorb}).

Using our systematic approach, we were able to conduct the strategy
on the model $\Gmod$, proving a conjecture of
Bostan, Bousquet-Mélou and Melczer in \cite{bostan2018counting}.
We detail the proof in Section~\ref{sec:meth-prov-algebr}
as an illustration of the strategy.
Moreover, for a family of models with large steps $(\calH_n)_n$ whose
orbits are conjectured to be finite, we were able to conjecture
that $X^{i+1} Y^{j+1}$ (which appears in place of $XY$ in equations
of the form (\ref{eq:eqfunc_g2}) when considering a starting point
$(i,j)$ other than $(0,0)$) admits a decoupling for several $(i,j)$
(namely, $(n-1,0)$ and $((n+1)k-1,k-1)$ for every $k$). We successfully proved
the algebraicity for several of these starting points for $n \le 4$, hinting
a possibly infinite family of algebraic models with arbitrarily large steps
(see Appendix~E of \cite{BHorb}).

\section{A Galois structure on the orbit} \label{sec:galo-struct-orbit}

The proofs and constructions in Sections \ref{sec:galo-struct-orbit},
\ref{sec:constr-invar-deco} and \ref{sec:meth-prov-algebr} are detailed in our
upcoming paper~\cite{BHorb}. We consider a weighted model $\calW$
with a non-univariate step polynomial. We denote by $k$ the field $\C(S(x,y))$.
Recall also that $\K = \overline{\C(x,y)}$, and that
if $M|L$ is a subextension of $\K|L$, a $L$-algebra automorphism
$\sigma\colon M \rightarrow M$ is a ring homomorphism such that
$\sigma_{|L} = \id_L$. We denote by $\Aut(M|L)$ the group
of $L$-algebra automorphisms of $M$.

We first endow the orbit with a group action as follows. If $\sigma : \K \rightarrow \K$ is a
$\C$-algebra automorphism, define its action on a pair $(u,v) \in \K \times \K$ by 
$\sigma \cdot (u,v) = (\sigma \cdot u, \sigma \cdot v)$.

\begin{lem}[Lemmas~3.7 and 3.8 in \cite{BHorb}]\label{lem:hom_graph}
  The orbit is stable under the action of $k(x)$ and $k(y)$-algebra automorphisms
  of $\K$, which all preserve the relations $\sim^x$ and $\sim^y$.
\end{lem}

This lemma has a field theoretic counterpart:
define $k(\cO)$ to be the subextension of $\K|k$ generated by the coordinates of the pairs
of $\cO$.
\begin{thm}[Theorem~3.9 in \cite{BHorb}]\label{thm:ext_gal}
  The field extensions $k(\cO) | k(x)$ and $k(\cO) | k(y)$ are Galois.
\end{thm}

We denote by $G_x=\Aut(k(\cO)|k(x))$ and $G_y=\Aut(k(\cO)|k(y))$ their respective Galois groups,
and by $G_{xy}$ their
intersection $G_x \cap G_y$ (which is the Galois group of the extension $k(\cO) | k(x,y)$).
We recall that the algebraic extension $k(\cO)|k(x)$ is \emph{Galois} if
$k(x)$ coincides with the subfield of $k(\cO)$ formed
by the elements fixed by every automorphism in $G_x$.

\begin{defi}[Group of the walk]
We define the \emph{group of the walk} $G = \left<G_x,G_y\right>$ to be the subgroup
of $k$-algebra automorphisms of $k(\cO)$ generated by $G_x$ and $G_y$.
\end{defi}

It is easy to see from its definition that $G$ acts by graph automorphisms on
the graph of $\cO$, and that its action is faithful.
Moreover, while the group $G$ is \emph{a priori} not finitely generated,
the left cosets $G_x/G_{x,y}$ and $G_y/G_{x,y}$ are of finite cardinal,
respectively $d_x = \deg_X \Ktld$ and $d_y = \deg_Y \Ktld$ (Lemma~3.14 in~\cite{BHorb}).
We then fix $I_x = \{\id, \iota^x_1, \dots, \iota^x_{d_x-1}\}$
and $I_y = \{\id, \iota^y_1, \dots, \iota^y_{d_y-1}\}$ two respective
sets of representatives for these two cosets.
\begin{thm}[Theorem~3.16 in \cite{BHorb}] \label{thm:grp_trans}
  The subgroup $\left<I_x,I_y\right>$ of $G$ acts transitively on $\cO$.
\end{thm}
Thus, the orbit $\cO$ is realized as the action of a finite set of automorphisms on
  the pair $(x,y)$, completing the analogy with the small steps setting.

\begin{exa}[Examples~3.10 and~3.18 in \cite{BHorb}]
  For $\calW$ a model with small steps that has both positive and
  negative steps in each direction, $k(\cO) = \C(x,y)$. Therefore, $G_{xy} = 1$,
  of index two in $G_x$ and $G_y$. Hence, $G_x = \left<\psi\right>$
  and $G_y = \left<\phi\right>$ with $\psi^2 = \phi^2 = 1$,
  and we find $G = \left<\phi,\psi\right>$. The identities
  $\psi(h(x,y)) = h(\Psi(x,y))$ and $\phi(h(x,y)) = h(\Phi(x,y))$
  yield an isomorphism between $G$ and the group of small steps
  (\S~\ref{sec:intro-group-orbit}).
\end{exa}

\begin{exa}[continuing Example~\ref{exa:orb_o12}]
  For $\Gmod$, we saw that $k(\cO) = \C(x,y,x_1)$, with $x_1$ algebraic of
  degree $2$ over $\C(x,y)$.
  Hence,  $G_{xy} = \left<\tau\right>$ with $\tau^2=1$,
  and after some computation we find $G_x = \left<\tau,\tau'\right> \simeq \Z / 2 \Z \times \Z / 2 \Z$
  and $G_y = \left<\tau,\sigma\right> \simeq S_3$ for $\tau'^2=1$
  and $\sigma^3=1$. In the end, $G = \left<\tau,\tau',\sigma\right>\simeq S_4$,
  which in this particular case coincides with the full group of graph automorphisms of $\cO$.
\end{exa}

\section{Construction of invariants and decoupling} \label{sec:constr-invar-deco}
\subsection{Fractions as elements of $k(\cO)$} \label{sec:fract-as-elem}

To apply the Galois framework of Section~\ref{sec:galo-struct-orbit} to
the construction of rational invariants and decoupling, we define an evaluation
of some fractions of $\C(X,Y,t)$ into $k(\cO)$.
  Its definition relies crucially on the fact that the kernel polynomial
  $\Ktld(X,Y.t)$ is irreducible in $\C[X,Y,t]$ (Lemma~3.10 in~\cite{BHorb}).

\begin{defi}
  We call a fraction $H(X,Y,t)$ of $\C(X,Y,t)$ \emph{regular} if
  the denominator of $H$ is not divisible by $\Ktld(X,Y,t)$.
\end{defi}
Note that fractions
of $\C(X,Y)$, $\C(X,t)$ and $\C(Y,t)$ are automatically regular because $\Ktld(X,Y,t)$ is
irreducible and trivariate by assumption on $\calW$.

\begin{defi}
  If $(u,v)$ is a pair of the orbit and $H(X,Y,t)$ is a regular fraction
  of $\C(X,Y,t)$, define its evaluation on $(u,v)$ to be
  $H_{(u,v)} = H(u,v,\invS) \in k(\cO)$.

  The evaluation on a pair of the orbit naturally extends to
  $\C$-linear combinations of pairs of the orbit (called \emph{$0$-chains}):
  \[
    \text{for } c = \sum_{(u,v) \in \cO} c_{u,v} (u,v), \quad \text{ define } \quad
    H_{c} = \sum_{(u,v) \in \cO} c_{u,v} H_{(u,v)}
  \]
\end{defi}

\begin{prop}[Proposition~3.23 in \cite{BHorb}]\label{prop:liftorb}
  The evaluation homomorphism sending a regular fraction
  $H$ to its evaluation $H_{(x,y)}$ maps bijectively $\C(X,t)$ to $k(x)$,
  $\C(Y,t)$ to $k(y)$ and $\C(X,Y)$ to $k(x,y)$. 
  A regular fraction evaluates to $0$ if and only if its numerator
  is divisible by $\Ktld(X,Y,t)$.
\end{prop}

Thus, we can consider regular fractions as elements of the field $k(\cO)$, so as to
benefit from our Galois-theoretic formalism. The homomorphism induces a relation on
regular fractions: two regular fractions of $\C(X,Y,t)$ are called \emph{Galois-equivalent}
if their evaluations induce the same element in $k(\cO)$.
Like the $t$-equivalence (Definition~\ref{defi:t_equiv}), this
  equivalence relation induces notions of invariants and decoupling:
\begin{itemize}
\item A pair of regular fractions $(I(X,t),J(Y,t))$ that are Galois equivalent is called a pair
  of \emph{Galois invariants}. By Proposition~\ref{prop:liftorb}, the evaluation
  homomorphism gives a correspondence between pairs of Galois invariants and
  elements of the subfield $k(x) \cap k(y)$ of $k(\cO)$, which we
  denote by $\kinv$, whose elements are fixed by $G$.
\item Likewise, we say that a regular fraction $H(X,Y,t)$ admits a \emph{Galois decoupling}
  pair $(F(X.t),G(Y,t))$ if $H$ is Galois-equivalent to $F + G$. As above,
  a fraction $H$ admitting a Galois decoupling corresponds through the
  evaluation homomorphism to a fraction $h$ in $k(x,y)$ that writes
  as $h = f + g$ for some $f \in k(x)$ (fixed by $G_x$) and $g \in k(y)$ (fixed by $G_y$).
\end{itemize}

Proposition~3.23 in \cite{BHorb} implies that $t$-equivalent regular
fractions are Galois equivalent. Therefore, the existence of rational
$t$-invariants or $t$-decoupling of a fraction $H$ is conditioned to the existence of
their Galois counterparts (of which we give a complete treatment in
the next two subsections).  Once we have obtained non-constant Galois
invariants or decoupling, we simply check if the Galois-equivalences
involved are also $t$-equivalences, so that we obtain non-constant
$t$-invariants and $t$-decoupling. If one of these two steps fails, then we
know that non-constant
rational $t$-invariants or $t$-decoupling of $H$ do not exist.

\subsection{Galois invariants} \label{sec:galois-invariants}

A pair of constant invariants $(F(t),F(t))$
is mapped by the evaluation homomorphism at $(x,y)$ to an element of $k$.
Therefore, the existence of non-constant Galois invariants is reduced
to the field-theoretic question of whether the inclusion $k \subset \kinv$ is proper or not.
This question is answered through the following result, which is a special
instance of a theorem proved in the more general context of finite
algebraic correspondences by Fried in~\cite{FriedPoncelet}, which
  we translate in the context of walks. This extends Theorem~4.6 in~\cite{BBMR16}
  and Corollary~4.6.11 in~\cite{FIM} to the large steps case.
\begin{thm}[Theorem~4.3 in \cite{BHorb}]
  The following statements are equivalent:
  \begin{enumerate}
  \item the orbit $\cO$ is finite,
  \item $G$ is a finite group,
  \item there exists a pair of non-constant Galois invariants.
  \end{enumerate}
\end{thm}

When finite, the orbit is described by two explicit polynomials that cancel the
left and right coordinates.
In this case, the extension $\kinv | k$ is purely transcendental
of transcendence degree $1$, and any nonconstant coefficient of these
polynomials is a generator, making the construction of rational invariants systematic
(e.g. Equation~(\ref{eq:inv_rat_gmod}) for $\Gmod$).
Finally, as $k(\cO)^G = \kinv$, then $k(\cO) | \kinv$ is a finite
Galois extension with Galois group $G$.

\subsection{Galois decoupling}\label{sec:decoupling}

We now assume that the orbit is finite. Given a regular fraction $H(X,Y,t)$,
we want to find a criterion for whether it admits a
Galois decoupling, and if it does, to compute it. To this end, we define a notion of
  decoupling in the orbit.

\begin{defi}[Definition~5.7 in \cite{BHorb}]\label{def:orb_decoupling}
  Let $(\gamma_x,\gamma_y,\alpha)$ be a tuple of $0$-chains such that
  $(x,y) = \gamma_x + \gamma_y + \alpha$ (with $(x,y)$ being considered
  as a vertex in the orbit). This is called a \emph{decoupling of $(x,y)$ in the orbit}
  if for every regular $H(X,Y,t)$ the following conditions hold:
  (1) $H_{\gamma_x} \in k(x)$, $H_{\gamma_y} \in k(y)$ and (2) $H_{\alpha} = 0$ when
  $H$ admits a Galois decoupling.
\end{defi}

In the proof of Theorem~4.11 in~\cite{BBMR16}, the authors construct an
explicit decoupling of
$(x,y)$ for the cyclic orbits of small steps models. We extend their
result to an arbitrary finite orbit using our Galois-theoretic framework
and graph homology.

\begin{thm}[Theorem~5.10 in \cite{BHorb}]\label{thm:decoupl_orbit}
  If the orbit is finite, the pair $(x,y)$ always admits a decoupling $(\gamma_x,\gamma_y,\alpha)$
  in the orbit (in the sense of Definition~\ref{def:orb_decoupling}).
\end{thm}

Thanks to this result, the question of the existence of the Galois decoupling of
a regular fraction may be decided through an evaluation:

\begin{prop} \label{prop:eval_dec}
  If $(\gamma_x,\gamma_y,\alpha)$ is a such a tuple, then
  a regular fraction $H(X,Y,t)$ admits a Galois decoupling if and only if $H_{\alpha} = 0$,
  and the decoupling is given by $H_{(x,y)} = H_{\gamma_x} + H_{\gamma_y}$.
\end{prop}

\begin{figure}[ht]
  \scalebox{0.5}{
    \input{o12_ldn.tikz}
  }
  \caption{The orbit $\cO_{12}$ of $\Gmod$ in two perspectives, illustrating
      a distance transitivity property (the $0$-chains $X_i$ and $Y_i$ are the
    sums of vertices in their respective regions).}
  \label{fig:orb_o12}
\end{figure}
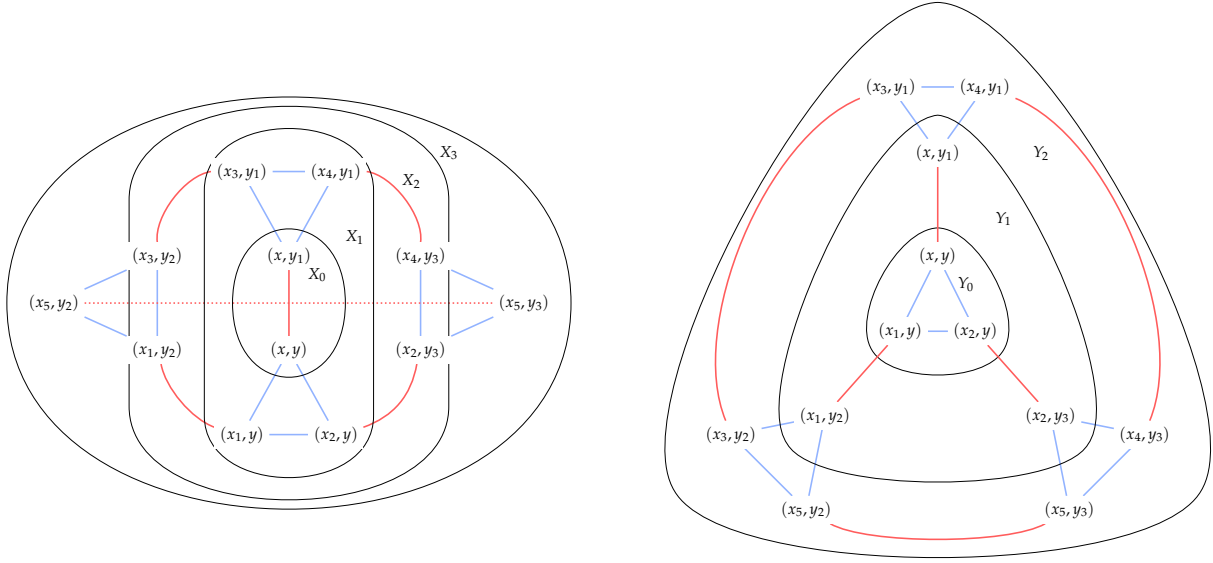

Note however that for an arbitrary $0$-chain $c$, it is not always convenient
to compute the evaluation $H_c$, because the coordinates of elements of
the orbit are random algebraic elements. A friendlier family for
computer algebra is composed of $0$-chains of the form
\[
  c = \sum_{(u,v) \in \cO, P(u)=0} (u,v) \quad \text{or} \quad c = \sum_{(u,v) \in \cO, P(v)=0} (u,v)
\]
with $P$ a polynomial over $\C(x,y)$, which we call \emph{symmetric chains}. It is easy
to compute the evaluation of regular fractions over symmetric chains via \emph{Newton's identities}.
Some symmetric chains are presented as \emph{level lines} for a well chosen
  distance in the graph of the orbit. They are denoted by $X_i$ and $Y_i$
in Figure~\ref{fig:orb_o12}.

Assuming a \emph{distance transitivity property} on the graph of the orbit,
we refine Theorem~\ref{thm:decoupl_orbit} by showing that the pair $(x,y)$ admits
a decoupling $(\gamma_x,\gamma_y,\alpha)$ with
$\gamma_x$ and $\gamma_y$ composed of level lines, whose expression is explicit
(Theorem~5.34 in~\cite{BHorb}).
This assumption was verified for all finite orbits arising from weighted models with
steps in $\{-1,0,1,2\}$ (Figure~10 in~\cite{bostan2018counting})
which includes the one of Figure~\ref{fig:orb_o12}. Other
families of orbits have been checked such as the ones arising from \emph{Hadamard}
(Section~11 of~\cite{bostan2018counting})
and \emph{Tandem} models (Section~3.2 of~\cite{BM2021plane}).

\begin{exa}\label{exa:decoupl_xy_o12}
  For the orbits of the same type as in Figure~\ref{fig:orb_o12},
  the pair $(x,y)$ admits a decoupling in terms of the symmetric chains $X_i$ and $Y_i$. It reads
  \begin{align*}
    (x,y) &= \left(\frac{X_0}{2}-\frac{X_1}{8}+\frac{X_2}{8}\right)+\left(\frac{Y_0}{4}-\frac{Y_1}{4}
            \right) + \alpha.
  \end{align*}
  We evaluated the fraction $XY$ on it to obtain its Galois decoupling~(\ref{eq:ex_decoupl_xy}).
\end{exa}

\section{An example: the model $\Gmod$ is algebraic} \label{sec:meth-prov-algebr}
We illustrate here with the model $\Gmod$ the generic strategy for proving
algebraicity of large steps models with small backward steps described in
Section~\ref{sec:invar-deco}.
First, recall the equation found for the generating function of walks on $\Gmod$ in
Example~\ref{ex:func_eq_gmod}:
\[
  \Ktld(X,Y,t) Q(X,Y,t) = XY - t (X^2 + 1) Q(X,0,t) - t Q(0,Y,t) + t Q(0,0,t). \qquad (\ref{eq:eqfunc_g2})
\]
Note that the generating function $Q(X,Y,t)$ has polynomial coefficients, hence
the left hand-side of the equation is $t$-equivalent to $0$,
so the functional equation translates into
\begin{equation} \label{eq:xy_equation}
  XY \equiv \left(t(X^2+1) Q(X,0,t) - t Q(0,0) \right) + t Q(0,Y,t).
\end{equation}
Moreover, using the decoupling of $(x,y)$ in the orbit of $\Gmod$ (Example~\ref{exa:decoupl_xy_o12}),
Proposition~\ref{prop:eval_dec} gives a Galois decoupling of the fraction $XY$, which is checked to
be the $t$-decoupling
\begin{equation} \label{eq:ex_decoupl_xy}
  XY \equiv - \frac{3 \lambda X^2 t   - \lambda t  - 4 X}{4 t (X^2 + 1) }
  + \frac{- \lambda Y - 4}{4 Y}.
\end{equation}
Combining Equations~(\ref{eq:xy_equation}) and (\ref{eq:ex_decoupl_xy}), we obtain the $t$-equivalence
\[
  \left(t(X^2+1) Q(X,0,t) - t Q(0,0,t) \right)
  +\frac{3 \lambda X^2 t   - \lambda t  - 4 X}{4 t (X^2 + 1) }
  \equiv \frac{- \lambda Y - 4}{4 Y}- t Q(0,Y,t),
\]
which gives a first pair of invariants $P_1=(I_1(X,t),J_1(Y,t))$.
Note that $J_1(Y,t)$ has a pole at $Y=0$, so this pair does not satisfy the conditions
of Lemma~\ref{lem:inv_lemma}.

The orbit of the model $\Gmod$ being finite, we also obtain automatically the following
pair $P_2=(I_2(X,t),J_2(Y,t))$ of Galois invariants, which we check to be $t$-invariants:
\begin{equation} \label{eq:inv_rat_gmod}
  \resizebox{.9\hsize}{!}{$
    \left(\frac{\left(-\lambda^{2} \,X^{3}-1  \,X^{4}-X^{6}+X^{2}+1\right) t^{2}-X^{2} \lambda  \left(X^{2}-1\right) t +X^{3}}{t^{2} X \left(X^{2}+1 \right)^{2}},
      \frac{-t \,Y^{4}+\lambda  t Y +Y^{3}+t}{Y^{2} t}\right).$
  }
\end{equation}

In order to find a pair of invariants satisfying the conditions of Lemma~\ref{lem:inv_lemma},
the heuristic is to combine the pairs of invariants $P_1$ and $P_2$ using
ring operations in order
to remove their poles both in $X$ and $Y$, by examining their Taylor
  expansions in their respective variables. Unlike the previous steps, the pole
  elimination is not systematic and requires a case by case treatment.
This leads us to define $P_3=(I_3(X,t),J_3(Y,t))$ to be
\[
  \resizebox{.9\hsize}{!}{$
    P_2\left(P_1-\frac{\lambda}{4}\right) -P_1^{3}+\left(2 t Q \! \left(0, 0\right)-\frac{\lambda}{4}\right) P_1^{2}+\left(2 t \frac{\partial Q}{\partial y}\! \left(0, 0\right)-t^{2} Q \! \left(0, 0\right)^{2}+\frac{5 \lambda^{2}}{16}\right) P_1$}.
\]
Using the functional equation, we are now able to check that this pair
of $t$-invariants indeed satisfies the conditions of
Lemma~\ref{lem:inv_lemma}. Therefore, there exists a power series
$A(t)$ in $\C[[t]]$ such that $I_3(X,t) = J_3(Y,t) = A(t)$. These are
equations with one catalytic variable for $Q(X,0,t)$ and $Q(0,Y,t)$
that satisfy the assumptions of Theorem~3 in~\cite{BMJ}. This allows
us to conclude that these series are algebraic over $\C(X,Y,t)$, so
that the same holds for the generating function $Q(X,Y,t)$ of the
model $\Gmod$. Following the method of~\cite{BMJ}, we found an
explicit minimal polynomial for the series $Q(0,0,t)$ of degree 32
with coefficients in $\Q(\lambda,t)$, proving in particular the
algebraicity conjecture on the excursion series of the two models
in~\cite{bostan2018counting} (lines~2 and~3 in Table~4, which are the
reversed models of $\calG_0$ and $\calG_1$ but sharing the same
excursion series).

\acknowledgements
We warmly thank Mireille Bousquet-Mélou for her invaluable advice and proofreading.

\printbibliography

\end{document}

%% file: article.tikzstyles
\tikzstyle{paire}=[fill=white, draw=none, shape=rectangle]
\tikzstyle{text_box}=[fill=white, draw=black, shape=rectangle]
\tikzstyle{decision}=[rectangle, minimum width=3cm, minimum height=1cm, text centered, draw=black]
\tikzstyle{arrow_flou}=[fill=none, draw={rgb,255: red,191; green,191; blue,191}, -Stealth, dashed]
\tikzstyle{clear_text}=[fill=white, draw=none, shape=rectangle, text={rgb,255: red,150;green,150; blue,150}]
\tikzstyle{start_point}=[fill=black, draw=black, shape=circle, radius=1pt]
\tikzstyle{ldn_text}=[fill=none, draw=black, shape=rectangle]
\tikzstyle{label}=[fill=white, draw=none, shape=circle]

\tikzstyle{simx}=[draw={rgb,255: red,146; green,179; blue,255}, ->, very thick]
\tikzstyle{simy}=[draw={rgb,255: red,255; green,92; blue,92}, ->, very thick]
\tikzstyle{usimx}=[draw={rgb,255: red,146; green,179; blue,255}, -, very thick]
\tikzstyle{usimy}=[draw={rgb,255: red,255; green,92; blue,92}, -, very thick]
\tikzstyle{light_grid}=[-, draw={rgb,255: red,210; green,210; blue,210}, thin, dashed]
\tikzstyle{walk_step}=[-, draw={rgb,255: red,162; green,43; blue,209}, -Stealth, very thick]
\tikzstyle{limits}=[-, pattern={Lines[angle=45,distance={5pt}]}, pattern color={rgb,255: red,180; green,180; blue,180}, draw=none]
\tikzstyle{arrow}=[thick, ->, >=Stealth]
\tikzstyle{fill_y_0}=[-, fill opacity=0.1, tikzit fill=none, fill={rgb,255: red,0; green,0; blue,255}]
\tikzstyle{fill_x_0}=[-, fill opacity=0.1, tikzit fill=none, fill={rgb,255: red,255; green,0; blue,0}]

\tikzset{edge from parent/.style=
{draw, edge from parent path={(\tikzparentnode.south)
-- +(0,-8pt)
-| (\tikzchildnode)}},
blank/.style={draw=none}}

%% file: g2_mod.tikz
\begin{tikzpicture}
	\begin{pgfonlayer}{nodelayer}
		\node [style=none] (0) at (0, 0) {};
		\node [style=none] (1) at (-2, -2) {};
		\node [style=none] (2) at (2, -2) {};
		\node [style=none] (3) at (2, 0) {};
		\node [style=none] (4) at (4, 2) {};
		\node [style=none] (5) at (0, 2) {};
		\node [style=none] (lamb) at (2, -0.5) {$\lambda$};
		\node [style=none] (6) at (0.5, 1.75) {};
		\node [style=none] (7) at (-1.25, -1.75) {};
	\end{pgfonlayer}
	\begin{pgfonlayer}{edgelayer}
	\draw[step=2.0,dashed,thin,color={rgb,255: red,210; green,210; blue,210}] (-2,-2) grid (4,2);
		\draw [style={axis_grid}] (0.center) to (5.center);
		\draw [style={axis_grid}] (0.center) to (4.center);
		\draw [style={axis_grid}] (0.center) to (3.center);
		\draw [style={axis_grid}] (0.center) to (2.center);
		\draw [style={axis_grid}] (0.center) to (1.center);
	\end{pgfonlayer}
\end{tikzpicture}

%% file: ex_walk_g2.tikz
\begin{tikzpicture}[scale=1]
	\begin{pgfonlayer}{nodelayer}
		\node [style=none] (0) at (-2, -2) {};
		\node [style=none] (1) at (-2, 8) {};
		\node [style=none] (2) at (0, 8) {};
		\node [style=none] (3) at (0, 0) {};
		\node [style=none] (4) at (8, 0) {};
		\node [style=none] (5) at (8, -2) {};
		\node [style=none] (6) at (4, 2) {};
		\node [style=none] (7) at (2, 4) {};
		\node [style=none] (8) at (0, 2) {};
		\node [style=none] (10) at (4, 6) {};
		\node [style=none] (11) at (4, 4) {};
		\node [style=none] (12) at (2, 2) {};
		\node [style=none] (13) at (4, 0) {};
		\node [style=none] (14) at (6, 0) {};
		\node [style=none] (15) at (4.75, 3) {};
		\node [style=none] (16) at (4.75, 5) {};
		\node [style=none] (17) at (2, 5.5) {};
		\node [style=none] (18) at (0.25, 3.75) {};
		\node [style=none] (19) at (1.25, 2.5) {$\lambda$};
		\node [style=none] (20) at (5, 0.5) {$\lambda$};
	\end{pgfonlayer}
	\begin{pgfonlayer}{edgelayer}
	\draw[step=2,color=gray,dashed] (0,0) grid (8,8);
		\draw [style=limits] (2.center)
			 to (3.center)
			 to (4.center)
			 to (5.center)
			 to (0.center)
			 to (1.center)
			 to cycle;
		\draw [style={walk_step}] (3.center) to (6.center);
		\draw [style={walk_step}] (6.center) to (11.center);
		\draw [style={walk_step}] (11.center) to (10.center);
		\draw [style={walk_step}] (10.center) to (7.center);
		\draw [style={walk_step}] (7.center) to (8.center);
		\draw [style={walk_step}] (8.center) to (12.center);
		\draw [style={walk_step}] (12.center) to (13.center);
		\draw [style={walk_step}] (13.center) to (14.center);
	\end{pgfonlayer}
\end{tikzpicture}

%% file: clas_small.tikz
\begin{tikzpicture}
\Tree [.{Models: 79} [.{$|G| < \infty$: 23} {Algebraic: 4} {D-finite: 19} ] [.{$|G| = \infty$ : 56} {D-algebraic : 9} {D-transcendental : 47} ] ]
\end{tikzpicture}

%% file: clas_large_steps.tikz
\begin{tikzpicture}
\Tree [.{Models: 13 110} [.{$|\cO| < \infty$: 13+227} {Algebraic: 2?} {D-finite: 227+4+7?} ] [.{$|\cO| = \infty$: 12 870} {non D-finite : $\ge$ 12 854 } ] ]
\end{tikzpicture}

%% file: o12_ldn.tikz
\begin{tikzpicture}[scale=2]
	\begin{pgfonlayer}{nodelayer}
		\node [style=paire] (0) at (9.5, 4) {$(x,y_1)$};
		\node [style=paire] (1) at (8.25, 5.75) {$(x_3,y_1)$};
		\node [style=paire] (2) at (10.75, 5.75) {$(x_4,y_1)$};
		\node [style=paire] (3) at (9.5, 1.25) {$(x,y)$};
		\node [style=paire] (4) at (8.5, -0.75) {$(x_1,y)$};
		\node [style=paire] (5) at (10.5, -0.75) {$(x_2,y)$};
		\node [style=paire] (6) at (6.5, -3) {$(x_1,y_2)$};
		\node [style=paire] (7) at (4, -3.5) {$(x_3,y_2)$};
		\node [style=paire] (8) at (6, -5.5) {$(x_5,y_2)$};
		\node [style=paire] (9) at (12.5, -3) {$(x_2,y_3)$};
		\node [style=paire] (10) at (13, -5.5) {$(x_5,y_3)$};
		\node [style=paire] (11) at (15, -3.5) {$(x_4,y_3)$};
		\node [style=paire] (12) at (-7.75, 1.25) {$(x,y_1)$};
		\node [style=paire] (13) at (-9, 3.5) {$(x_3,y_1)$};
		\node [style=paire] (14) at (-6.5, 3.5) {$(x_4,y_1)$};
		\node [style=paire] (15) at (-7.75, -1.25) {$(x,y)$};
		\node [style=paire] (16) at (-9, -3.5) {$(x_1,y)$};
		\node [style=paire] (17) at (-6.5, -3.5) {$(x_2,y)$};
		\node [style=paire] (18) at (-11.25, -1.25) {$(x_1,y_2)$};
		\node [style=paire] (19) at (-11.25, 1.25) {$(x_3,y_2)$};
		\node [style=paire] (20) at (-14, 0) {$(x_5,y_2)$};
		\node [style=paire] (21) at (-4.25, -1.25) {$(x_2,y_3)$};
		\node [style=paire] (22) at (-1.5, 0) {$(x_5,y_3)$};
		\node [style=paire] (23) at (-4.25, 1.25) {$(x_4,y_3)$};
		\node [style=none] (24) at (9.5, 2) {};
		\node [style=none] (25) at (7.75, -1.25) {};
		\node [style=none] (26) at (11.25, -1.25) {};
		\node [style=none] (27) at (9.5, 5) {};
		\node [style=none] (28) at (5.5, -3.75) {};
		\node [style=none] (29) at (13.5, -3.75) {};
		\node [style=none] (39) at (-0.25, 0) {};
		\node [style=none] (40) at (-15.25, 0) {};
		\node [style=none] (41) at (10.25, 0.5) {$Y_0$};
		\node [style=none] (42) at (11.25, 2.25) {$Y_1$};
		\node [style=none] (43) at (12.25, 4) {$Y_2$};
		\node [style=none] (44) at (-7, 0.75) {$X_0$};
		\node [style=none] (45) at (-6, 1.75) {$X_1$};
		\node [style=none] (46) at (-4.5, 3.25) {$X_2$};
		\node [style=none] (47) at (-3.5, 4) {$X_3$};
		\node [style=none] (48) at (9.5, 8) {};
		\node [style=none] (49) at (2.5, -5) {};
		\node [style=none] (50) at (16.5, -5) {};
		\node [style=none] (51) at (-6.25, 0) {};
		\node [style=none] (52) at (-9.25, 0) {};
		\node [style=none] (53) at (-10, 3) {};
		\node [style=none] (54) at (-5.5, 3) {};
		\node [style=none] (55) at (-5.5, -3) {};
		\node [style=none] (56) at (-10, -3) {};
		\node [style=none] (57) at (-12, 2.75) {};
		\node [style=none] (58) at (-3.5, 2.75) {};
		\node [style=none] (59) at (-3.5, -2.75) {};
		\node [style=none] (60) at (-12, -2.75) {};
	\end{pgfonlayer}
	\begin{pgfonlayer}{edgelayer}
		\draw [style={none}] (26.center)
			 to [bend left=60, looseness=0.75] (25.center)
			 to [bend left=60, looseness=0.75] (24.center)
			 to [bend left=60, looseness=0.75] cycle;
		\draw [style=usimx] (4) to (3);
		\draw [style=usimx] (3) to (5);
		\draw [style=usimx] (5) to (4);
		\draw [style=usimx] (1) to (0);
		\draw [style=usimx] (0) to (2);
		\draw [style=usimx] (2) to (1);
		\draw [style=usimx] (6) to (7);
		\draw [style=usimx] (7) to (8);
		\draw [style=usimx] (8) to (6);
		\draw [style=usimx] (9) to (10);
		\draw [style=usimx] (10) to (11);
		\draw [style=usimx] (11) to (9);
		\draw [style=usimy] (4) to (6);
		\draw [style=usimy] (5) to (9);
		\draw [style=usimy, bend right, looseness=0.50] (8) to (10);
		\draw [style=usimy] (3) to (0);
		\draw [style=usimy, bend right=45, looseness=0.75] (1) to (7);
		\draw [style=usimy, bend left=45, looseness=0.75] (2) to (11);
		\draw [style=usimx] (16) to (15);
		\draw [style=usimx] (15) to (17);
		\draw [style=usimx] (17) to (16);
		\draw [style=usimx] (13) to (12);
		\draw [style=usimx] (12) to (14);
		\draw [style=usimx] (14) to (13);
		\draw [style=usimx] (18) to (19);
		\draw [style=usimx] (19) to (20);
		\draw [style=usimx] (20) to (18);
		\draw [style=usimx] (21) to (22);
		\draw [style=usimx] (22) to (23);
		\draw [style=usimx] (23) to (21);
		\draw [style=usimy, bend left, looseness=0.75] (16) to (18);
		\draw [style=usimy, bend right] (17) to (21);
		\draw [style=usimy, dotted] (20) to (22);
		\draw [style=usimy] (15) to (12);
		\draw [style=usimy, bend right=45, looseness=0.75] (13) to (19);
		\draw [style=usimy, bend left=45, looseness=0.75] (14) to (23);
		\draw [style={none}] (29.center)
			 to [bend right=60, looseness=0.50] (27.center)
			 to [bend right=60, looseness=0.50] (28.center)
			 to [bend right=60, looseness=0.50] cycle;
		\draw [style={none}] (39.center)
			 to [bend left=90, looseness=1.25] (40.center)
			 to [bend left=90, looseness=1.25] cycle;
		\draw [style={none}] (50.center)
			 to [bend right=60, looseness=0.50] (48.center)
			 to [bend right=60, looseness=0.50] (49.center)
			 to [bend right=60, looseness=0.50] cycle;
		\draw [style={none}] (51.center)
			 to [bend left=90, looseness=2.25] (52.center)
			 to [bend left=90, looseness=2.25] cycle;
		\draw [style={none}] (53.center)
			 to [bend left=90, looseness=1.25] (54.center)
			 to (55.center)
			 to [bend left=90, looseness=1.25] (56.center)
			 to cycle;
		\draw [style={none}] (57.center)
			 to [bend left=90] (58.center)
			 to (59.center)
			 to [bend left=90] (60.center)
			 to cycle;
	\end{pgfonlayer}
\end{tikzpicture}